\documentclass[10pt]{amsart}
\usepackage{amssymb,a4wide}
\usepackage{amsthm}
\usepackage{graphicx}
\usepackage[T2A]{fontenc}
\usepackage[latin1]{inputenc}
\usepackage{tensor}
\usepackage{amscd}
\usepackage{enumerate}

\newcommand{\weg}[1]{}
\newcommand{\diag}{\mathrm{diag}}

\newcommand{\Id}{\mathrm{Id}}

\theoremstyle{plain}
\newtheorem{thm}{Theorem}
\newtheorem*{thm*}{Theorem}
\newtheorem{lem}{Lemma}

\theoremstyle{definition}

\theoremstyle{remark}
\newtheorem{rem}{Remark}

\begin{document}

\title{On the number of nontrivial projective transformations of closed  manifolds}
\author{Vladimir S. Matveev}
\address{Friedrich-Schiller-Universit\"at Jena, 
07737 Jena Germany  \ {\bf Email:}  vladimir.matveev@uni-jena.de}
\date{}

\begin{abstract} 
We show that for a closed Riemannian manifold
the quotient of the group of projective transformations  by the  group of isometries contains at most two elements unless the metric has constant positive sectional curvature or every projective transformation is an affine transformation. 

\end{abstract} 
 
\maketitle

 \vspace{-0.3cm}

\hfill {\small \it  Dedicated to Anatoly Timofeevich Fomenko  on his seventieth birthday.}

 \vspace{0.3cm}

\vspace{1ex} 
Let $(M, g)$ be a smooth closed Riemannian manifold of dimension $n\ge 2$. By \emph{projective transformation} of $(M, g)$ we understand a diffeomorphism $\phi:M\to M$ that sends geodesics viewed as unparameterized curves to geodesics. Projective transformations of $(M, g)$ naturally form a group which we denote by $\operatorname{Proj}$. The group of isometries of $(M, g)$,  which we denote by $\operatorname{Iso}$,  forms a subgroup of $\operatorname{Proj}$. The following theorem, which is the main result of our paper,    answers the following natural question: how big 
can be $\operatorname{Proj}/\operatorname{Iso}$?

\begin{thm} \label{thm} Let $(M, g)$ be a smooth connected closed Riemanian 
manifold  of dimension $n\ge 2$. 
Suppose   $\left|\operatorname{Proj}/\operatorname{Iso}\right|> 2$. Then, the 
sectional curvature of $g$  is  positive constant or 
 every projective transformation is an affine transformation, i.e., preserves the Levi-Civita  connection of $g$. 
 \end{thm}   

In other words, if a connected closed Riemannian manifold of dimension $\ge 2$ whose sectional curvature is not positive constant admits two nonaffine  projective transformations $\phi$ and $\psi$ then the  diffeomorphisms   $\phi\circ \psi, \psi\circ \phi,  \phi^{-1} \circ  \psi, \phi\circ \psi^{-1}$ are    isometries. 

Both possibilities in the conclusion of Theorem \ref{thm} can happen. It is well known, see for example \cite[Example 2]{gonone}, that for the standard sphere 
$(S^n, g_{standard})$ we have 
$$
\operatorname{Proj}/\operatorname{Iso}= \operatorname{SL}(n+1, \mathbb{R})/\operatorname{SO}(n+1, \mathbb{R})
$$
so the set  $\operatorname{Proj}/\operatorname{Iso}$ contains infinitely many elements. 

Note though that certain quotients of  odd-dimensional spheres of constant positive sectional curvature have $\left|\operatorname{Proj}/\operatorname{Iso}\right|=1.$ This follows from \cite[Theorem 1]{topology}, and a concrete 3-dimensional  example that can be generalized for all odd dimensions  is in \cite[\S 1.4]{mounoud}.  

Let us also recall   an example (\cite[Example 4]{gonone})  such that  $\left|\operatorname{Proj}/\operatorname{Iso}\right|=\infty,$ but each projective transformation is actually an affine transformation. Consider the standard torus $T^2= \mathbb{R}^2/\mathbb{Z}^2$, where the action of the group  $\mathbb{Z}^2$ is generated by the standard translations $(x,y)\mapsto (x +1, y)$ and $(x,y)\mapsto (x, y+1)$ along the standard basis vectors.  The standard flat metric on $\mathbb{R}^2$ induces a 
metric on $T^2$ which we denote by  $g$.  
Consider now the standard action of $SL(2, \mathbb{Z})$  on $\mathbb{R}^2$. It induces a faithful  action 
of $SL(2, \mathbb{Z})$ on $T^2$ which evidently preserves the Levi-Civita  connection of $g$. 
Hence, $\left|\operatorname{Proj}/\operatorname{Iso}\right|=\infty.$ Note that though the example is two-dimensional and flat, it is easy to extend it to nonflat manifolds of higher  dimensions by taking  direct products with compact manifolds.

 \vspace{1ex}

 Let us now construct an example of a closed Riemannian manifold of arbitrary dimension $n\ge 2$ such that $\operatorname{Proj}/\operatorname{Iso}$ contains two elements 
 and such that it does not admit affine nonisometric transformations. Two-dimensional version of this example is in \cite[\S 1.3]{difgeo} and in \cite[Example 4]{gonone}.   We consider the direct product  $S^1 \times M \times S^1$, where $S^1$ is the circle and  $M$ is an arbitrary closed connected manifold of dimension $n-2$. We denote by $x$ and $z$ the standard cyclic coordinates on the first and the second $S^1$; we assume that $x,y\in (\mathbb R \mod 1)$. We will denote by $y_1,...,y_{n-2}$ local coordinates on the manifold $M$.  We arbitrary choose    a Riemannian metric $g$ on $M$ and a smooth nonconstant 
 1-periodic function  $f:\mathbb{R} \to \mathbb{R}$ such that $f>1$. 
 
 Consider the following metric on $S^1 \times M\times  S^1$: 
 
 $$
\left( f(x)- \tfrac{1}{f(z)}\right)\left(f(x)- 1\right)  dx^2 + (f(x)-1)\left(1- \tfrac{1}{f(z)}\right) \sum g_{ij} dy^idy^j + \left( f(x)- \tfrac{1}{f(z)}\right) \left(1- \tfrac{1}{f(z)}\right)  dz^2.
 $$
 
Next, consider the diffeomorphism $\phi:S^1\times M \times S^1 $ given by $\phi(x,y,z)= (z,y,x)$  (where $y$ denotes   a point on $M$).  This pullback with respect to this diffeomorphism is given by 

 $$
\left( f(x)- \tfrac{1}{f(z)}\right)\left(f(x)- 1\right) \tfrac{f(z)}{ f(x)^2 } dx^2 + (f(x)-1)\left(1- \tfrac{1}{f(z)}\right) \tfrac{f(z)}{f(x)} \sum g_{ij} dy^idy^j + \left( f(x)- \tfrac{1}{f(z)}\right) \left(1- \tfrac{1}{f(z)}\right) \tfrac{f(z)^2}{f(x)} dz^2, 
 $$
which is projectively equivalent (i.e., has the same geodesics) to the above metric by the Levi-Civita theorem \cite{LC}, and the diffeomorphism $\phi$ is therefore a projective transformation.

\begin{rem} The diffeomorphism $\phi:(x,y,z)\mapsto (z,y,x)$ used in the example above is not orientable. One can easily modify the example such that the projective transformation is orientable. Indeed, 
 in dimensions $\ge 3$, one can   choose  $(M,g)$ such that it  admits an orientation-reversing isometry $\alpha:M\to M$, 
and instead of diffeomorphism $\phi$ consider the diffeomorphism $(x,y,z)\mapsto (z,\alpha(y),x)$. Essentially the same idea works in dimension 2 as well: one takes the function $f$  such that it is even and superpose the diffeomorphism $\phi$ with an orienation-reversing isometry $(x,y,z)\mapsto (-x,y,z)$. 
\end{rem} 

The following special case of Theorem \ref{thm} 
 is  due  to A. Zeghib \cite{zeghib}, where also some  previous results in this direction are listed, see also the introductions to \cite{difgeo} and to \cite{gonone}  for an overview and for the  history of the problem.  Zeghib \cite[Theorem 1.3]{zeghib}  has proved   Theorem \ref{thm}  under a stronger assumption $\left|\operatorname{Proj}/\operatorname{Iso}\right|>2n.$ Actually, our proof follows the  lines  of 
 the  Zeghib's proof  and is based on his results and ideas; we will clearly explain the additional argument that allowed us to improve his result.

\begin{rem}   Theorem \ref{thm} remains correct if one replaces closeness (of our manifold $(M,g)$) by completeness.  
 This improved statement is  based on certain  nontrivial 
 ideas and calculations that are invented in  \cite{CEMN} and will be published after or in that paper. 
\end{rem}

\vspace{1ex} 
\noindent {\bf Proof of Theorem \ref{thm}.} Within the proof we 
 assume that $(M,g)$ is a closed connected Riemannian manifold of dimension at least  $2$  admitting at least one projective transformation which is not an affine transformation. We also assume that the sectional curvature is not positive constant. Our goal is to show $\left|\operatorname{Proj}/\operatorname{Iso}\right|\le 2.$

Consider the metrization equation from \cite[Theorem 2.2]{EM}. The precise formula  of this equation is not important for us and its introducing requires work we do not want to invest, we refer to \cite{EM} for details. We will list here the properties of this equation and its solutions which will be used in the proof. 

\begin{enumerate}[(I)] \item The metrization equation is a (homogeneous) linear system of PDE, 
 so   its solution space which we denote by $\operatorname{Sol}$ is a linear vector space. By \cite[Theorem 2]{degree} (or alternatively  \cite[Thereom 16]{difgeo}, \cite[Theorem 1]{KioMat2010}, \cite[Corollary 5.2]{mounoud}; the two-dimensional version follows from  \cite{kiyohara}),
  under our assumptions, $\dim \operatorname{Sol}\le 2$. 

\item In a local coordinate system, 
the solutions of the metrization equation can be viewed by  matrices\footnote{As geometric objects they are weighted symmetric $(2,0)$-tensors; in particular their  pullback is well-defined.} whose components  are functions of  the coordinates. 
  Nondegenerate  (i.e., with nowhere vanishing determinant)  
    solutions  $\sigma^{ij}$
  correspond, in a local coordinate system,  to  metrics (of arbitrary signature) 
   projectively equivalent to $g$ (i.e., having the same geodesics with $g$). The correspondence is given  by the formula
  \begin{equation}\label{L}
  \sigma^{ij}= g^{-1} \, |\det g|^{\tfrac{1}{n+1}} \ \textrm{ and } \  g^{-1} := \det|\sigma|\  \sigma^{ij}.
  \end{equation}
From the formula we clearly see that positively definite $\sigma^{ij}$ correspond to positively definite, i.e., Riemannian,   metrics.

\item Metrization equations are projectively invariant, so for any projective transformation  $\phi$ 
  the pullback $\phi^*\sigma$ of a solution is a solution.

   \item If a solution $\sigma$ is nondegenerate at every point, then $\sigma^{-1} \bar \sigma:= (\sigma^{-1})_{is}\bar \sigma^{sj}$ 
   (where $\bar \sigma$ is also a solution), is a well-defined (1,1)-tensor field. If $\sigma$ and $\bar \sigma$ correspond to the  metrics $g$ and $\bar g$ by the formula \eqref{L},  we have $\sigma^{-1} \bar \sigma= \left|\tfrac{\det  \bar g}{\det  g}\right|^{\tfrac{1}{n+1}} \bar g^{js}g_{si}.$
   
\end{enumerate}

Suppose $\phi$ is a projective stransformations. Take a basis $\sigma, \bar \sigma $ in $\operatorname{Sol}$ and consider the pullbacks $\phi^*\sigma$, $\phi^*\bar \sigma$. They also belong to $\operatorname{Sol}$ and are therefore linear combinations of the basis solutions  $\sigma$ and $\bar \sigma$; we denote the coefficients as below: 
$$
\begin{pmatrix} \phi^*\sigma \\ \phi^*\bar \sigma\end{pmatrix} = \begin{pmatrix} a & b \\ c & d \end{pmatrix} \begin{pmatrix} \sigma\\ \bar \sigma\end{pmatrix} = \begin{pmatrix} a\sigma  &+ & b\bar \sigma  \\ c \sigma  &+& d\bar \sigma \end{pmatrix}.
$$ 
We denote the matrix $\begin{pmatrix} a & b \\ c & d \end{pmatrix}$ above by $A$ or by $A_\phi$.  
The  mapping  from $\operatorname{Proj}$ to $GL(2, \mathbb{R})$ given by $\phi\mapsto A_\phi$ 
is actually a representation and if $A_\phi$ is  the identity matrix  then $\phi$ is an  isometry.    The composition  $\psi\circ \phi$ of two projective transformations corresponds to the product of matrices $A_\psi $ and $A_\phi$ is the reverse order:
$$
\psi\circ \phi \mapsto A_\phi A_\psi. 
$$

In \cite{zeghib} it was shown that, in our assumptions, if $A_\phi$ has  real eigenvalues, then 
      $\phi$ is  an isometry. The result is nontrivial and of huge importance for us since further we may assume that the matrices $A_\phi$ have nonreal eigenvalues. 
 
 \vspace{1ex}

Let us now consider the case when  the matrix   $A=A_\phi$ has complex-conjugate nonreal  eigenvalues. 
Depending on the sign of the determinant,  by the choice of the basis in $\operatorname{Sol}$,  we may therefore assume that the matrix $A$ is as in one of the following two cases: 
$$
A= C \begin{pmatrix} \cos \alpha & \sin \alpha \\ - \sin \alpha & \cos \alpha  \end{pmatrix}  \ \ \textrm{or} \ \ 
A= C \begin{pmatrix} \cos \alpha & \sin \alpha \\  \sin \alpha & -\cos \alpha  \end{pmatrix} 
$$ 
with $C>0$.   The first case is when $\det A>0$ and the second is when $\det A<0$. Moreover, without loss of generality we may assume in the first case that the basis $\sigma, \bar \sigma$
 is such that the metric $g$ corresponds to the basis solution $\sigma$.  

Let us show   that the first case is impossible, unless $C=1$ and $\alpha$ is such that $\begin{pmatrix} \cos \alpha & \sin \alpha \\ - \sin \alpha & \cos \alpha  \end{pmatrix}= \begin{pmatrix} 1 & 0\\ 0& 1\end{pmatrix}.$
This is precisely the argument   overseen by Zeghib in his paper. 

In order to do this, let us take an arbitrary point $p\in M$ and a  basis in $T_pM$ such that 
$$
\sigma= \diag(1,...,1) \ \ \textrm{and} \ \ \bar\sigma= \diag(s_1,..., s_n). 
$$
The existence of such a basis is trivial since $g$ is positively definite.   

Next, consider $$\phi^* \sigma, \ \phi\circ \phi^*\sigma=  \phi^*(\phi^*(\sigma)), \ \phi\circ \phi\circ \phi^* \sigma = \phi^*(\phi^*(\phi^*(\sigma))),...,\underbrace{\phi\circ ... \circ \phi^*}_{\textrm{$k$ times}}\sigma ,... \ . $$ Since the matrix corresponding to the  superposition $\underbrace{\phi\circ ... \circ \phi}_{\textrm{$k$ times}}$ is  $\begin{pmatrix} \cos \alpha & \sin \alpha\\ -\sin \alpha  & \cos \alpha\end{pmatrix}^k=  \begin{pmatrix} \cos k\alpha & \sin k\alpha\\ -\sin k\alpha  & \cos k\alpha\end{pmatrix}$, we have  that at our point $p$  the matrix $\sigma^{-1} \left(\underbrace{\phi\circ ... \circ \phi^*}_{\textrm{$k$ times}}\sigma\right)$ is

\begin{equation} \label{eigen}
\sigma^{-1}(C^k(\cos k\alpha \sigma + \sin k\alpha \bar \sigma)     ) = C^k\diag\bigl( \cos k\alpha + s_1  \sin k\alpha,...,\cos k\alpha + s_n  \sin k\alpha\bigr).
\end{equation}
We will need the following simple  lemma:

\begin{lem} \label{lem} 
 Suppose  for all  $k\in \mathbb{N}$ we have $\cos k\alpha + s  \sin k\alpha>0$. Then,  
 $\alpha$ is such that $$\begin{pmatrix} \cos \alpha & \sin \alpha \\ - \sin \alpha & \cos \alpha  \end{pmatrix} = \begin{pmatrix} 1 & 0\\ 0& 1\end{pmatrix}.$$  
\end{lem}  

\vspace{1ex} 
\noindent {\bf Proof of  Lemma \ref{lem}.} First observe  that $\cos k\alpha + s  \sin k\alpha>0$ is the 1st coordintate  of the  $k\alpha$-rotation of the vector $\begin{pmatrix} 1\\ s\end{pmatrix} $   around the origin. 
If $\begin{pmatrix} \cos \alpha & \sin \alpha \\ - \sin \alpha & \cos \alpha  \end{pmatrix} \ne \begin{pmatrix} 1 & 0\\ 0& 1\end{pmatrix}$, the set of the points of $\mathbb{R}^2$ of the form 
 $\begin{pmatrix} \cos k\alpha & \sin k\alpha \\ - \sin k\alpha & \cos k\alpha  \end{pmatrix} \begin{pmatrix} 1\\ s\end{pmatrix}$ is either dense on a circle centered at the origin, or coincides with the set of 
  vertices of a  regular (may be degenerate, i.e., containing only two vertices)
  polygon containing the origin. In both case there exists a point such that its first coordinate is nonpositive and we get a contradiction.   {\bf Lemma \ref{lem} is proved.}

\vspace{1ex}   
We now continue the proof of Theorem \ref{thm}. Combining  Lemma \ref{lem} with \eqref{eigen}, we see that if  $\begin{pmatrix} \cos \alpha & \sin \alpha \\ - \sin \alpha & \cos \alpha  \end{pmatrix} \ne  \begin{pmatrix} 1 & 0\\ 0& 1\end{pmatrix}$, there exists   $k$ such that the (diagonal) 
 matrix of  $\sigma^{-1} (\phi \circ ...\circ \phi^*\sigma)$ has   a nonpositive eigenvalue. But since $\sigma$ corresponds to the metric and is therefore positively definite at all points of the manifold which implies of course that its pullback is also positively definite, 
 this is impossible. Thus, $A_\phi= C \Id$. Since closed  manifolds do not  admit nontrivial  homotheties, $C=1$ and $\phi$ is an isometry. 
 
Finally, we obtained that for every projective transformation $\phi$ such that it is not an isometry the matrix   $A_\phi$ has two complex conjugate nontrivial eigenvalues    and negative determinant. Then, for two such (nonisometric) projective transformations $\phi$ and $\psi$ there superpositions $\phi\circ \psi $ is an isometry, 
since the product of two matrices $A_\psi$ and $A_\phi$ with negative determinants 
 has positive determinant. Then, all projective transformations such that they are not isometries lie in the same 
 equivalence class of $\operatorname{Proj}/\operatorname{Iso}$ which implies that the number of elements in  $\operatorname{Proj}/\operatorname{Iso}$ is at most 2. {\bf Theorem \ref{thm} is proved. }


\begin{thebibliography}{99}

\bibitem{CEMN}
D. Calderbank, M. Eastwood, V. S,  Matveev, K.  Neusser,   \emph{C-projective geometry: background and  open problems, } in preparation. 
   

  
\bibitem{EM} 
M. Eastwood,  V. S. Matveev,  \emph{ Metric connections in projective differential geometry,}
 Symmetries and Overdetermined Systems of Partial Differential Equations (Minneapolis, MN, 2006), 339--351,
 IMA Vol. Math. Appl.,    {\bf
144}(2007),   Springer, New York.  




\bibitem{KioMat2010} V. Kiosak, V. S. Matveev, \emph{Proof Of The
    Projective Lichnerowicz Conjecture For Pseudo-Riemannian Metrics
    With Degree Of Mobility Greater Than Two}, Comm. Mat. Phys. {\bf  297}(2010),
  no. 2, 401--426.  

\bibitem{kiyohara} K. Kiyohara, \emph{Compact Liouville surfaces}, J. Math. Soc. Japan {\bf 43}(1991), 555--591.


 \bibitem{LC}
 T. Levi-Civita, {\it Sulle trasformazioni delle equazioni
 dinamiche}, Ann. di Mat., serie $2^a$, {\bf 24}(1896), 255--300.


\bibitem{gonone} V. S. Matveev, \emph{Beltrami problem, Lichnerowicz-Obata conjecture and applications of integrable systems in differential geometry}, Tr. Semin. Vektorn. Tenzorn. Anal, {\bf 26}(2005), 214--238. 

\bibitem{degree} V. S. Matveev, \emph{ On degree of mobility of complete metrics}, Compt. Math., {\bf 43}(2006), 221--250.

 

\bibitem{topology} V. S. Matveev,
\emph{ Three-dimensional manifolds having metrics with the same geodesics,}
    Topology     {\bf 42}(2003) no. 6,  1371--1395.

 \bibitem{difgeo} V. S. Matveev, \emph{  Proof of Projective Lichnerowicz-Obata Conjecture}, J. of Differential Geometry, {\bf 75}(2007), 459--502.
    
    \bibitem{mounoud} V. S. Matveev, P. Mounoud,
  \emph{Gallot-Tanno Theorem for closed incomplete pseudo-Riemannian
    manifolds and applications},  Ann. Glob. Anal. Geom. {\bf 38}(2010), 259--271.

\bibitem{zeghib} A. Zeghib, \emph{
On discrete projective transformation groups of Riemannian manifolds,} arXiv:1304.6812v1.


\end{thebibliography}
\end{document}